\documentclass[10pt,a4paper]{article}
\usepackage[utf8]{inputenc}
\usepackage[english]{babel}

\usepackage[left=2cm,right=2cm,top=2cm,bottom=2cm]{geometry}

\usepackage{cite}
\usepackage{amsmath,amssymb,amsfonts}
\usepackage{graphicx}
\usepackage{textcomp}

\usepackage{multirow}
\usepackage[electronic]{ifsym}
\usepackage{indentfirst,latexsym}
\usepackage{bm}
\usepackage{bbm}
\usepackage{mathrsfs}
\usepackage{amsmath,amsfonts,amscd,amssymb}
\usepackage{pifont}
\usepackage[T1]{fontenc}
\usepackage{mathbbol}

\usepackage{subfigure}

\usepackage{xcolor}

\newcommand{\blue}[1]{{\color[RGB]{0,0,255} #1}}

\usepackage{url}

\DeclareMathOperator{\dif}{d}         

\DeclareMathOperator{\MOD}{mod}
\renewcommand{\mod}{\MOD}

\renewcommand{\vec}[1]{\bm{#1}}

\newcommand{\mfloor}[1]{ \left\lfloor {#1} \right\rfloor }
\newcommand{\mceil}[1]{ \left\lceil {#1} \right\rceil }
\newcommand{\mpair}[2]{ \left\langle {#1}, {#2} \right\rangle}

\newcommand{\set}[1]{\left\{ #1 \right\}}
\newcommand{\seq}[1]{\langle #1 \rangle}
\newcommand{\abs}[1]{\left| #1 \right|}
\newcommand{\braket}[2]{ \langle #1 | #2 \rangle}

\newcommand{\Zernpoly}[2]{\operatorname{Z}_{#1}^{#2}}
\newcommand{\Radipoly}[2]{\operatorname{R}_{#1}^{#2}}

\newcommand{\scrd}[2]{{#1}_{\mathrm{#2}}}

\newcommand{\BigO}[1]{\mathscr{O}\left({#1}\right)}

\usepackage{algorithm}         
\usepackage{algorithmicx}
\usepackage{algpseudocode}
\usepackage{xcolor}

\newcommand{\Algr}{\textbf{Algorithm}~}
\newcommand{\Fig}{FIGURE~}
\newcommand{\Tab}{TABLE~}

\algnewcommand\algorithmicswitch{\textbf{switch}}
\algnewcommand\algorithmiccase{\textbf{case}}
\algnewcommand\algorithmicdefault{\textbf{default}}
\algnewcommand\algorithmicassert{\texttt{assert}}
\algnewcommand\Assert[1]{\State \algorithmicassert(#1)}%
\algdef{SE}[SWITCH]{Switch}{EndSwitch}[1]{\algorithmicswitch\ #1\ \algorithmicdo}{\algorithmicend\ \algorithmicswitch}%
\algdef{SE}[CASE]{Case}{EndCase}[1]{\algorithmiccase\ #1}{\algorithmicend\ \algorithmiccase}%
\algdef{SE}[DEFAULT]{Default}{EndDefault}[1]{\algorithmicdefault\ #1}{\algorithmicend\ \algorithmicdefault}%

\makeatletter
\renewcommand{\ALG@name}{Algorithm}
\newenvironment{breakablealgorithm}
{
	\begin{center}
		\refstepcounter{algorithm}
		\setlength{\baselineskip}{15pt} 
		\renewcommand{\caption}[2][\relax]{
			\hrule height.9pt depth0pt \kern3pt
			{\raggedright\textbf{\ALG@name~\thealgorithm} ##2\par}%
			\ifx\relax##1\relax 
			\addcontentsline{loa}{algorithm}{
				\protect\numberline{\thealgorithm}##2}%
			\else 
			\addcontentsline{loa}{algorithm}{
				\protect\numberline{\thealgorithm}##1}%
			\fi
			\kern2pt\hrule\kern2pt
		}
	}{
		\kern3pt\hrule\relax
	\end{center}
}

\usepackage{caption}
\usepackage{listings}
\newcommand{\cpvar}[1]{\texttt{#1}}

\newcommand{\ProcName}[1]{\textsc{#1}}

\DeclareMathOperator{\Tc}{\mathrm{T}}   
\DeclareMathOperator{\Tf}{\mathrm{TF}}  
\DeclareMathOperator{\spacomp}{\mathcal{S}} 
\DeclareMathOperator{\timcomp}{\mathcal{T}} 
\newcommand{\tcc}[2]{\timcomp_{#1}{\left(#2\right)}}
\newcommand{\scc}[2]{\spacomp_{#1}{\left(#2\right)}}

\author{Hong-Yan Zhang\thanks{Corresponding author: Hong-Yan Zhang, email: hongyan@hainnu.edu.cn}, Yu Zhou and Zhi-Qiang Feng\\
School of Information Science and Technology, Hainan Normal University, Haikou 571158, China}
\title{\textbf{Balanced Binary Tree Schemes for Computing \\ Zernike Radial Polynomials}}
\date{Aug. 19, 2024}

\begin{document}
\maketitle

\begin{abstract}
\textit{Zernike radial polynomials} (ZRP) play a significant role in application areas such as optics design, imaging systems, and image processing systems. Currently, there are two kinds of numerical schemes for computing the ZRP automatically with computer programs: one is based on the definition in which the factorial operations may lead to the overflow problem and the high order derivatives are troublesome, and the other is based on recursion which is either unstable or with high computational complexity. In this paper, our emphasis is focused on exploring the \textit{balanced binary tree} (BBT) schemes for computing the ZRP: firstly an elegant formulae for computation is established; secondly  the recursive  and iterative algorithms based-on BBT are proposed; thirdly  the computational complexity of the algorithms are  analyzed rigorously; finally the performance of BBT schemes by testing the running time is verified and validated. Theoretical analysis shows that the computational complexity of \textit{balanced binary tree recursive algorithm} (BBRTA) and iterative algorithm are exponential and quadratic respectively, which coincides with the running time test very well.  Experiments show that the time consumption is about $1\sim 10$ microseconds with different computation platforms for the \textit{balanced binary tree iterative algorithm} (BBTIA), which is stable and efficient for real-time applications. In the sense of STEM education, the connection of the BBT and ZRP exhibits the beauty and applications of discrete mathematical structure behind the engineering problem, which is worthy of introducing to the college students, computer programmers and optics engineers.
\end{abstract}

\section{Introduction}

The Zernike radial polynomials (ZRP),  named after Frits Zernike, are important for atmospheric turbulence analysis\cite{Noll-1976}, aberration analysis in imaging system and optics design \cite{ZHY2023ZCP,Mahajan-1981,Diaz-2013,OptikShopTest-2007,Mathar-2008,Buhren2018,Berger2022-ZernikeAberr,BornWolf-1999,Jasssen-2010,Shakibaei-2013,Zemax-2008}, and image processing \cite{Chong-2003}. Mathematically, ZRP are a sequence of orthogonal polynomials which are defined on the unit disk and derived from the pupils of imaging system. Generally, for the radial degree  $n\in\mathbb{Z}^+ = \set{0, 1, 2, \cdots}$ and azimuthal order $m\in \mathbb{Z}$ such that $\abs{m}\le n$ and $n-m$ is even, the ZRP with double indices $\mpair{n}{m}$ are denoted by \cite{BornWolf-1999,Jasssen-2010,Shakibaei-2013}
\begin{equation} \label{eq-Rnm-def}
\begin{split}
 \Radipoly{n}{m}(\rho)
 &= \frac{1}{(\frac{n-\abs{m}}{2})!\rho^{\abs{m}}} \left[\frac{\dif}{\dif (\rho^2)} \right]^{\frac{n-\abs{m}}{2}}
    \left[ (\rho^2)^{\frac{n+\abs{m}}{2}} (\rho^2 -1)^{\frac{n-\abs{m}}{2}} \right] \\
&=\sum_{s=0}^{k} c_s \rho^{n-2s} = c_0\rho^n + c_1\rho^{n-2} + \cdots + c_{k} \rho^{n-2k}
\end{split}
\end{equation} 
in which the parameter for the number of terms is 
\begin{equation}\label{eq-k-def}
k = \frac{n-\abs{m}}{2}
\end{equation}
and the coefficients $c_s$ can be expressed by
\begin{equation} \label{eq-cs}
c_s = (-1)^s \binom{k}{s}\binom{n-s}{k}, \quad 0\le s \le k
\end{equation}
where
\begin{equation}
\binom{\alpha}{i} = \frac{\alpha(\alpha-1)\cdots(\alpha-i+1)}{i!}
=\prod^{i-1}_{t=0} \frac{\alpha-t}{i-t}
\end{equation}
is the binomial coefficient.
The symmetry characterized by  
\begin{equation} \label{eq-Rnm-symmetry}
\Radipoly{n}{m}(\rho) = \Radipoly{n}{-m}(\rho) = \Radipoly{n}{\abs{m}}(\rho)
\end{equation}
implies that it is sufficient to consider the non-negative $m$ for computing the value of $\Radipoly{n}{m}(\rho)$. It is easy to find that it is not wise for us to compute $\Radipoly{n}{m}(\rho)$ directly based on the definition with the equation \eqref{eq-Rnm-def} since the factorial operations may lead to the overflow problem and the high order derivatives are troublesome in practical computations and/or implementations with computer programs.

In the past decades of exploring the indirect computation methods for $\Radipoly{n}{m}(\rho)$, there are several recursive methods to compute $\Radipoly{n}{m}(\rho)$. In 1976, Kintner \cite{Kintner-1976} proposed the \textit{$n$-recursive} formulae 
\begin{equation} \label{eq-Rnm-Kintner}
\begin{aligned}
\Radipoly{n}{m}(\rho) 
&= \frac{1}{k_1}\left[(k_2\rho^2 + k_3)\Radipoly{n-2}{m}(\rho)
+ k_4\Radipoly{n-4}{m}(\rho)\right], \quad n\ge 4
\end{aligned}
\end{equation} 
where
\begin{equation}
\left\{
\begin{split}
k_1 &= K_1(n,m) = (n+m)(n-m)(n-2)/2, \\
k_2 &= K_2(n,m) = 2n(n-1)(n-2), \\
k_3 &= K_3(n,m) = -m^2(n-1)-n(n-1)(n-2),\\
k_4 &= K_4(n,m) = -n(n+m-2)(n-m-2)/2.\\
\end{split}
\right.
\end{equation}
However, the stopping condition is unknown at that time and the formula is singular when $k_1 = 0$. 
In 1989, Prata and Rusch  \cite{Prata-1989} proposed the following recursive scheme
\begin{equation} \label{eq-R-PR-v1}
\Radipoly{n}{m}(\rho) = \rho L_1 \Radipoly{n-1}{m-1}(\rho) + L_2\Radipoly{n-2}{m}(\rho), \quad n\ge 2
\end{equation}
with the coefficients
\begin{equation}
L_1 = \frac{2n}{m+n}, \quad L_2 = \frac{m-n}{m+n} = 1- L_1.
\end{equation}
and stopping condition
\begin{equation} \label{eq-R-sc-PR-v1}
\Radipoly{n}{m}(\rho)= \left\{
\begin{array}{ll}
\rho,  & n = 1; \\
1,     & n = 0.
\end{array}
\right.
\end{equation}
There are two significant points for the stopping condition:
\begin{itemize}
\item it has significant impacts on the  time complexity and space complexity of the computation;
\item the choice of stopping condition is not unique, which leads to different performance for the numeric schemes adopted for the objective of computation.
\end{itemize}
For computing the $\Radipoly{n}{m}(\rho)$, Chong et. al. proposed a simple stopping condition for recursive process in 2003 \cite{Chong-2003}
\begin{equation} \label{eq-R-sc}
\Radipoly{n}{m}(\rho)
=
\left\{
\begin{array}{ll}
\rho^m, & n = m; \\
\rho^m[(m+2)\rho^2 - (m+1)], & n = m + 2.
\end{array}
\right.
\end{equation}
In 2013, Shakibaei and Paramesran \cite{Shakibaei-2013} reformulated the recursive relation in \eqref{eq-R-PR-v1} by 
\begin{equation} \label{eq-R-PR}
\Radipoly{n}{m}(\rho) = \rho L_1\Radipoly{n-1}{\abs{m-1}}(\rho) + L_2\Radipoly{n-2}{m}(\rho)
\end{equation}
and derived an alternative recursive scheme
\begin{equation} \label{eq-R-SP}
\Radipoly{n}{m}(\rho) = \rho\left[\Radipoly{n-1}{\abs{m-1}}(\rho) + \Radipoly{n-1}{m+1}(\rho) \right] - \Radipoly{n-2}{m}(\rho), \quad n\ge 2
\end{equation} 
with the stopping condition
\begin{equation} \label{eq-R-sc-SP}
\Radipoly{n}{m}(\rho) =
\left\{
\begin{array}{ll}
 0, & n <m;\\
 \rho, & n = 1;\\
1, & n = 0; 
\end{array}
\right.
\end{equation}
via the properties of Chebyshev polynomials of the second kind. 
However, the computation process is rather slow with this recursive method. Chong et al. \cite{Chong-2003} proposed the following \textit{$m$-recursive scheme}
\begin{equation} \label{eq-Rnm-chong}
\begin{aligned}
\Radipoly{n}{m-4}(\rho) 
&=  \left(h_2 + \frac{h_3}{\rho^2} \right)\Radipoly{n}{m-2}(\rho) 
+ h_1\Radipoly{n}{m}(\rho), \quad m\ge 4
\end{aligned}
\end{equation}
in which $h_1, h_2$ and $h_3$ are functions of $n$ and $m$.
This $m$-recursive scheme is more efficient than the other recursive schemes for computing $\Radipoly{n}{m}(\rho)$.  However, $\rho = 0$ is a singular point in \eqref{eq-Rnm-chong} although $\Radipoly{n}{m}(\rho)$ is regular for all $\rho\in[0, 1]$. Thus the computation will be unstable if $\rho$ is small enough. 

In computer science, we know that the essence of recursion lies in two facts: there must be a stopping condition for the recursive procedure/function which calls itself; the recursion depth should not be large otherwise the computational complexity will be too large  due to the massive memory consumption and long time consumption caused by the push-pop stacking processes. 
For the available recursive schemes of computing $\Radipoly{n}{m}(\rho)$  at present, for large $n$ and difference of $n-m$, the computation complexity of these recursive schemes is rather high which limits their applications. Generally, for the recursive problem with single integer as argument, it is easy to convert the recursive formulae to a more efficient iterative counterpart. However, for the $\Radipoly{n}{m}(\rho)$ with two integers $n$ and $m$ as arguments, there is a lack of feasible method to convert the recursive formula to iterative versions. Although Kintner's $n$-recursive formula \eqref{eq-Rnm-Kintner} can be reformulated as an iterative formula, it is limited for $n\ge 4$; Chong's $m$-recursive scheme \eqref{eq-Rnm-chong} can be converted to its iterative version, however the singular point $\rho = 0$ will still exist. For the coupled recursive formulae \eqref{eq-R-PR-v1} (or \eqref{eq-R-PR} equivalently) and \eqref{eq-R-SP}, their  iterative implementations are still to be explored.

For the computational complexity of the recursive algorithms available,  
Shakibaei and Paramesran \cite{Shakibaei-2013} considered the time complexity by counting the number of addition and multiplication operations. However, their conclusion is arguable for some reasons: the space complexity of computation is ignored, the running time is not tested and the complexity just depends on the radial index $n$ instead of the double indices $n$ and $m$. It should be noted that the space complexity of recursive algorithm is usually exponential or more higher. In consequence, the recursive schemes are not suitable for real-time applications.  

In this paper, our objective is to explore an elegant formula to compute $\Radipoly{n}{m}(\rho)$ stably and propose novel recursive and iterative schemes with the help of \textit{balanced binary tree} (BBT) structure. Our main contributions lie in the following   perspectives:
\begin{itemize}
\item A novel formula for computing the ZRP is proposed, which stimulates the recursive and iterative schemes for the numerical computation.
\item The BBT structure of the novel formula is discovered, which deepens the understanding of the ZRP.
\item Both the recursive and iterative algorithms for computing  the ZRP are designed, which refreshes the state-of-the-art of the computational complexity.
\end{itemize}

The rest of this paper is organized as follows: Section \ref{sec-pre} deals with the preliminaries for developing objective algorithms; Section \ref{sec-recursive} discusses the BBT recursive scheme; Section \ref{sec-iterative} concerns the BBT iterative scheme; Section \ref{sec-v-and-v} focuses on the verification and validation; finally Section \ref{sec-summary} gives the conclusion for our work.

\section{Preliminaries} \label{sec-pre}

\subsection{A Special Kind of Difference Equation}

For the discrete difference equation
\begin{equation} \label{eq-Gp-de}
G_p  = \alpha G_{p-2} + \beta, \quad p = 2, 4, 6, \cdots
\end{equation}
with initial value $G_0$ where $\alpha, \beta$ are constants and $\alpha \neq 1$, its solution is
\begin{equation} \label{eq-ans-Gp-de}
G_p = \alpha^{\frac{p}{2}-1}G_0 + \frac{\beta}{\alpha -1}\left(\alpha^{\frac{p}{2}-1} -1 \right), \quad p = 2, 4, 6, \cdots
\end{equation}
This formula will be used in analyzing the computational complexity of computing $\Radipoly{n}{m}(\rho)$.

\subsection{Computing Powers of Real Number with Squaring}

For any $n\in \mathbb{Z}^+$ and real number $x\in \mathbb{R}$, the power $x^n$ can be computed fast by squaring, which is based on the following recursive formulae
\begin{equation}
x^n = \left\{
\begin{array}{ll}
x^{n-1}\cdot x, & n\ge 1, 2\nmid n;\\
x^{\frac{n}{2}}\cdot x^{\frac{n}{2}}, & n\ge 1, 2\mid n;\\
1, & n = 0.
\end{array}
\right.
\end{equation}
With this formulae, for $n\ge 1$ the power $x^n$ can be computed fast with the time complexity of $O(\log_2 n)$ since only $\mceil{\log_2 n}$ times of multiplicative operation is required. However, for $n=0$, $x^n$ is always $1$, thus the complexity will be $\BigO{1}$. The notations $\mathscr{O}(\cdot)$, $\scc{\ProcName{ProcName}}{\cdot}$,
$\tcc{+}{\cdot}$, $\tcc{*}{\cdot}$,
$\tcc{\ProcName{Alg}}{\cdot}$,  $\Tf(\cdot)$, $\Tc_{\ProcName{Alg}}(\cdot)$,
$\Tf_{\ProcName{Alg}}(\cdot)$, $\scc{\ProcName{Alg}}{\cdot}$
and so on about the computational complexity are introduced in the appendix, please see Appendix \ref{app-subsec-cc} for more details.

The iterative algorithm for computing $x^n$ is shown in \Algr \ref{alg-power-squaring}.

\begin{breakablealgorithm}
\caption{Compute the power $x^n$ with squaring method in an iterative way} \label{alg-power-squaring}
\begin{algorithmic}[1]
\Require variable $x\in \mathbb{R}$, variable $n\in \mathbb{Z}^+$
\Ensure the value of power $x^n$
\Function{CalcPower}{$x$, $n$}
\State $\cpvar{prod} \gets 1$;
\While{$n\ge 1$}
\If{$2 \nmid n$}
   \State $\cpvar{prod}\gets \cpvar{prod} \cdot x$;
\EndIf
\State $x \gets x\cdot x$;
\State $ n \gets n/2 $; 
\EndWhile
\State \Return \cpvar{prod};
\EndFunction
\end{algorithmic}
\end{breakablealgorithm}
Obviously, the time complexity of \Algr \ref{alg-power-squaring} is
\begin{equation} \label{eq-tcc-power}
\tcc{\ProcName{CalcPower}}{n} = 
\left\{
\begin{array}{ll}
\BigO{1}, & n = 0;\\
\BigO{\log_2 n}, & n\ge 1.
\end{array}
\right.
\end{equation}

For the purpose of computing the generalZRP $\Radipoly{n}{m}(\rho)$, it is necessary to investigate the functions $\Radipoly{m}{m}(\rho)$ and $\Radipoly{m+2}{m}(\rho)$ since they can be used to generate formulae for $\Radipoly{n}{m}(\rho)$ as stopping/initial conditions for recursive/iterative processes. In the recursive tree, there are two types of leaf nodes: one is type A specified by $\Radipoly{m}{m}(\rho)$, another is type B specified by $\Radipoly{m+2}{m}(\rho)$. 
\Algr \ref{alg-leaf-type-A} and \Algr \ref{alg-leaf-type-B} are used to compute the leaf nodes of type A with $\Radipoly{m}{m}(\rho)$ and type B with $\Radipoly{m+2}{m}(\rho)$ respectively.

\begin{breakablealgorithm}
\caption{Compute the value of a leaf node of type A with $\Radipoly{m}{m}(\rho) =\rho^{m}$ where $m\in \mathbb{Z}^+$.}
\label{alg-leaf-type-A}
\begin{algorithmic}[1]
\Require Radius $\rho \in [0,1]$, integer $m\in \mathbb{Z}^+$
\Ensure  The value of $\Radipoly{m}{m}(\rho)=\rho^m$.
\Function{CalcLeafNodeTypeA}{$\rho,m$}
\State \Return $\ProcName{CalcPower}(\rho,m)$;
\EndFunction
\end{algorithmic}
\end{breakablealgorithm}

\begin{breakablealgorithm}
\caption{Compute the value of a leaf node of type B with $\Radipoly{m+2}{m}(\rho) =(m+2) \rho^{m+2} - (m+1)\rho^m$ where $m\in \mathbb{Z}^+$.}
\label{alg-leaf-type-B}
\begin{algorithmic}[1]
\Require Radius $\rho \in [0,1]$, integer $m\in \mathbb{Z}^+$
\Ensure  The value of $\Radipoly{m}{m+2}(\rho)=\rho^m[(m+2)\rho^2 -(m+1)]$.
\Function{CalcLeafNodeTypeB}{$\rho,m$}
\State \Return $\ProcName{CalcPower}(\rho,m)\cdot((m+2)\rho^2 - (m+1))$;
\EndFunction
\end{algorithmic}
\end{breakablealgorithm}

The time flops for the leaf nodes of type A is given by
\begin{equation} \label{eq-tcc-leaf-type-A}
\begin{split}
\tcc{}{\ProcName{CalcLeafNodeTypeA}}
=\tcc{*}{\ProcName{CalcPower}(\rho,m)} 
= \left\{
\begin{array}{ll}
0, & m = 0;\\
\mceil{\log_2 m}, & m\ge 1
\end{array}
\right.
\end{split}
\end{equation}
As a comparison, the time flops for the leaf nodes of type B can be expressed by
\begin{equation} \label{eq-tcc-leaf-type-B}
\begin{split}
\tcc{}{\ProcName{CalcLeafNodeTypeB}}
=&\tcc{}{\ProcName{CalcPower}(\rho,m)} +  1 
+ \tcc{}{(m+2)\rho^2-(m+1)}  \\
=&
\left\{
\begin{array}{ll}
3, & m = 0;\\
3 + \mceil{\log_2 m}, & m\ge 1.
\end{array}
\right.
\end{split}
\end{equation}
Therefore, the computational complexity for the leaf nodes of type A and type B is always $\BigO{\log_2 m}$ for $m\ge 1$ or $\BigO{1}$ for $m=0$ when computing $\Radipoly{m}{m}(\rho)$ and $\Radipoly{m+2}{m}(\rho)$ with index $m$ and radius $\rho$ (or $\Radipoly{n}{n-2}(\rho)$ with index $n$ and radius $\rho$).

\section{Recursive Scheme for Radial Polynomials} \label{sec-recursive}

\subsection{Computing the Radial Polynomials Recursively}

Our novel recursive formulae is a combination of the recursive schemes in \cite{Shakibaei-2013} and \cite{Prata-1989}. The trick of the exploring is to reduce the difference of the up-down scripts appearing on the right hand side in \eqref{eq-R-PR} and \eqref{eq-R-SP}  with a common constant so as to get a balanced result, see \Tab 
\ref{tab-up-down-scripts}. The larger the difference of the up-down scripts in $\Radipoly{n}{m}(\rho)$ is, viz. $n-m$, the faster the recursive process is. For the two terms $\rho L_1\Radipoly{n-1}{\abs{m-1}}(\rho)$ and $L_2 \Radipoly{n-2}{m}(\rho)$  in the right hand side of \eqref{eq-R-PR}, the reductions of the difference, namely $(n'-m') -(n-m)$, are $0$ and $2$ where $\Radipoly{n'}{m'}(\rho)$ is in the right hand side of recursive formula. Obviously, the reductions are not equal. Similarly, for the three terms in the right hand side of 
\eqref{eq-R-PR}, the reductions of the difference are $0$, 
$2$ and $2$ respectively. These quantities are also not balanced.
It should be noted that if the reduction of difference is $0$, then  it is slow for the recursive process to satisfy the stopping condition.

\begin{table*}[htp]
\centering
\caption{Reduction of $n-m$ for the up-down scripts of $\Radipoly{n}{m}(\rho)$ in recursive relations}
\label{tab-up-down-scripts}
\begin{tabular}{crcclcc}
\hline
Eq.  & Expr. of $\Radipoly{n'}{m'}(\rho)$   & $m'$ & $n'$ &  $n'-m'$ & Reduction & Balanced Reduction?\\
\hline \hline
 \multirow{2}*{\eqref{eq-R-PR}}
& $\rho L_1\Radipoly{n-1}{\abs{m-1}}(\rho)$ &  $\abs{m-1}$ & $n-1$ & $n-m$  & $0$ ($m\ge 1$) & \multirow{2}*{No}\\
 & $L_2\Radipoly{n-2}{m}(\rho)$ &   $m$ & $n-2$ & $n-m-2$ & $2$ & \\
\hline
\multirow{3}*{\eqref{eq-R-SP}} & $\rho \Radipoly{n-1}{\abs{m-1}}(\rho)$ &  $\abs{m-1}$ & $n-1$ & $n-m$  & $0$ ($m\ge 1$) &  \multirow{3}*{No}\\
 & $\rho \Radipoly{n-1}{m+1}(\rho)$ & $m+1$ & $n-1$ & $n-m-2$  & $2$\\
 &  $-\Radipoly{n-2}{m}(\rho)$  & $m$ & $n-2$ & $n-m-2$ & $2$\\
\hline
\multirow{2}*{\eqref{eq-Rnm-recu}} & $\rho F_1 \Radipoly{n-1}{m+1}(\rho)$ 
& $m+1$ & $n-1$ & $n-m-2$ & $2$ & \multirow{2}*{Yes}\\
 & $F_2 \Radipoly{n-2}{m}(\rho)$ 
& $m$ & $n-2$ & $n-m-2$ & $2$ &  \\
\hline
\end{tabular}
\end{table*}

\begin{figure*}[htb]
\centering
\includegraphics[width=1.0\textwidth]{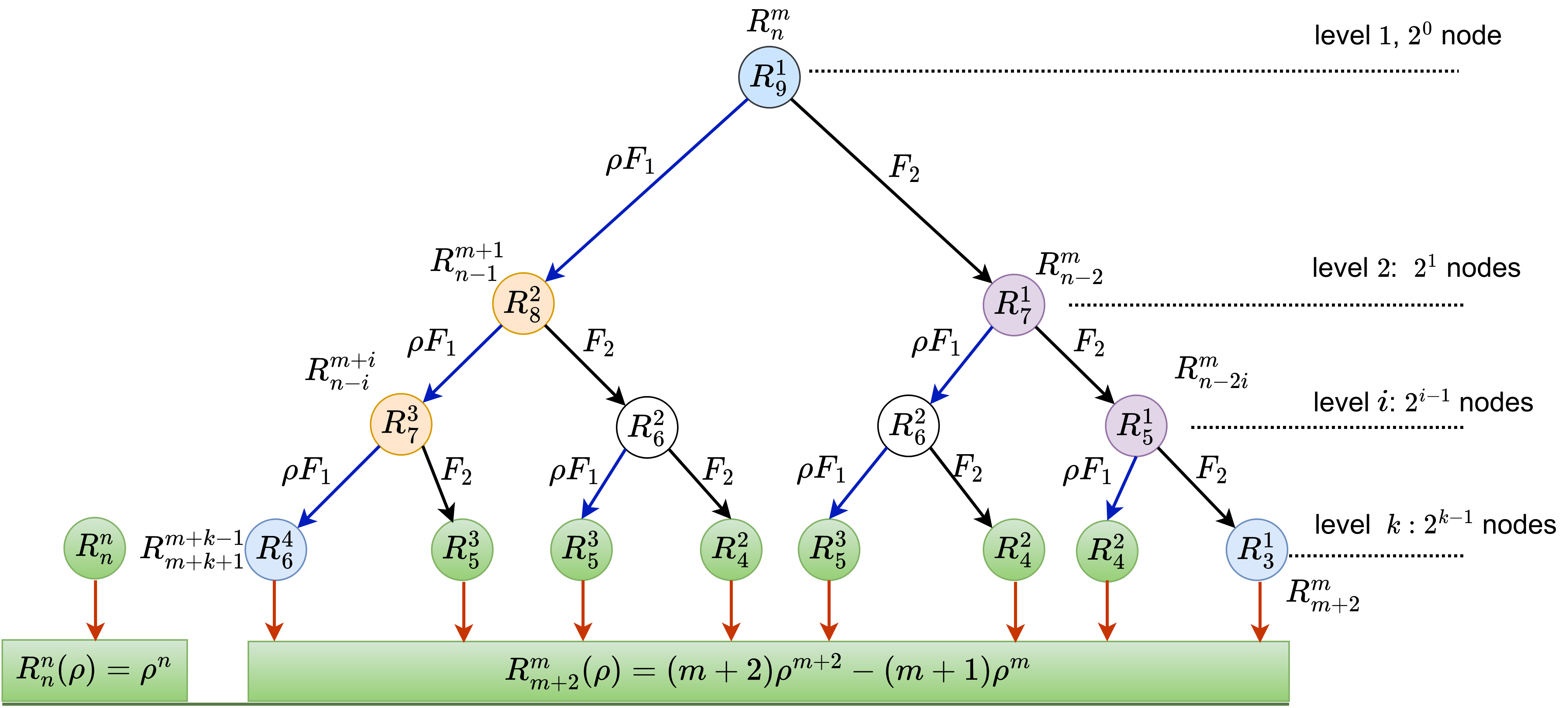} 
\caption{Principle of Recursive scheme for computing $\Radipoly{n}{m}(\rho)$ with BBT when $n\neq m\ge 0$}
\label{fig-recursive}
\end{figure*}

Let 
\begin{equation} \label{eq-F-nm}
\left\{
\begin{split}
F_1 &= F_1(n,m)=\frac{2n}{n-m}, \\
F_2 &=F_2(n,m)= -\frac{n+m}{n-m}= 1-F_1,
\end{split}
\right. \quad n\neq m
\end{equation}
by multiplying $L_1$ with \eqref{eq-R-SP} and eliminating the term $\rho L_1\Radipoly{n-1}{\abs{m-1}}(\rho)$ in \eqref{eq-R-PR}, we immediately have
\begin{equation} \label{eq-Rnm-recu}
\Radipoly{n}{m}(\rho) 
= \rho F_1 \Radipoly{n-1}{m+1}(\rho) + F_2 \Radipoly{n-2}{m}(\rho), \quad n-m \ge 2, m\ge 0.
\end{equation}
Evidently, the differences of up-down scripts in $\Radipoly{n-1}{m+1}(\rho)$ and $\Radipoly{n-2}{m}(\rho)$ are the same, i.e., 
\begin{equation}
(n-1)-(m+1) = (n-2) -(m) = n-m -2.
\end{equation} 
In this way, we can reduce the difference $n-m$ to $0$ or $2$ step by step with a common difference $d = 2$, and then the stopping condition \eqref{eq-R-sc} can be applied naturally. 
In conclusion, we  can obtain the following complete recursive scheme 
\begin{equation} \label{eq-R-BBT}
\left\{
\begin{array}{ll}
\Radipoly{n}{m}(\rho) = \Radipoly{n}{\abs{m}}(\rho), & m <0; \\
\Radipoly{n}{m}(\rho) = \rho F_1\Radipoly{n-1}{m+1}(\rho) + F_2\Radipoly{n-2}{m}(\rho),  &m \ge 0, n-m \ge 4;\\
\Radipoly{n}{m}(\rho) = \rho^m, & m \ge 0, n-m = 0;\\
\Radipoly{n}{m}(\rho) = \rho^m[(m+2)\rho^2-(m+1)], & m \ge 0, n-m = 2;\\
\end{array}
\right.
\end{equation}
For practical computation, we can always set $m$ with $\abs{m}$ at first. Before proceeding further, we give some interpretations about this recursive scheme:
\begin{itemize}
\item the recursive process corresponds to a BBT as shown in \Fig \ref{fig-recursive};
\item the stopping conditions consists of two kinds of leaves: 
     \begin{itemize}
     \item type A, which corresponds to $\Radipoly{m}{m}(\rho)$ for $n=m$;
     \item type B, which corresponds to $\Radipoly{m+2}{m}(\rho)$ for $n= m+2$;
     \end{itemize}
\item the root node corresponds to the objective polynomial $\Radipoly{n}{m}(\rho)$; 
\item computing the values of the leaf nodes with type A and type B  can be easily and fast performed with the stopping condition \eqref{eq-R-sc};
\item for $n\neq m\ge 0$, the computation of $\Radipoly{n}{m}(\rho)$ corresponds to traversing the BBT and  computing the value for each node in the tree according to \eqref{eq-F-nm} and \eqref{eq-Rnm-recu} dynamically;
\item the number of levels of the nodes in the binary tree is equal to the parameter $k$ in \eqref{eq-k-def};
\item the number of nodes in the $i$-th level is $2^{i-1}$ for $i\in \set{1, 2, \cdots, k}$ and the total number of nodes is $\displaystyle \sum^{k}_{i=1} 2^{i-1} =2^k -1$. 
\end{itemize}

\subsection{Balanced Binary Tree Recursive Algorithm (BBTRA)}

For the given recursive formulae and stopping condition, generally it is easy to design a recursive algorithm. 
With the help of the preceding two procedures  \ProcName{CalcLeafNodetypeA} and \ProcName{CalcLeafNodetypeB} specified by the stopping condition \eqref{eq-R-sc}, we can
design the \textit{balanced binary tree recursive algorithm} (BBTRA) as  listed in \Algr \ref{alg-Rnm-recur}.
\begin{breakablealgorithm}
\caption{Compute the radial function $\Radipoly{n}{m}(\rho)$  recursively}
\label{alg-Rnm-recur}
\begin{algorithmic}[1]
\Require Born-Wolf double indices $\mpair{n}{m}$, variable $\rho\in[0, 1]$
\Ensure the value of $\Radipoly{n}{m}(\rho)$
\Function{CalcRadiPolyBBTRA}{$n$, $m$, $\rho$}
\State $m \gets \abs{m}$;  \quad // $\Radipoly{n}{m}(\rho) = \Radipoly{n}{-m}(\rho)$. 
\State $p\gets n-m$;\quad // $p = 0, 2, 4, 6, \cdots$
\Switch{$p$}
\Case{$0$}  \quad // for the leaf node $\Radipoly{m}{m}(\rho) = \rho^m$.
 \State \Return $\ProcName{CalcLeafNodeTypeA}(\rho,m)$; 
\EndCase
\Case{$2$} \quad // for the leaf nodes $\Radipoly{m+2}{m}(\rho)$.
\State  \Return $\ProcName{CalcLeafNodeTypeB}(\rho, m)$; 
\EndCase
\Default \quad // for the node $\Radipoly{n}{m}=\Radipoly{m+p}{m}(\rho),p \ge 4$
   \State $F_1\gets 2n/(n-m), F_2 \gets 1- F_1$;
   \State \Return  $\rho\cdot F_1 \cdot \ProcName{CalcRadiPolyBBTRA}(n-1, m+1, \rho)+F_2 \cdot\ProcName{CalcRadiPolyBBTRA}(n-2, m, \rho)$;
\EndDefault
\EndSwitch
\EndFunction
\end{algorithmic}
\end{breakablealgorithm}

It should be noted that $p=n-m$ is the difference odd integers or even integers, which implies that $p=0, 2, 4, \cdots$. Therefore, there is no \fbox{\textbf{case 1}} in \Algr  \ref{alg-Rnm-recur}.

\subsection{Computational Complexity of BBTRA}

For the non-negative double indices $n$ and $m$ such that $n - m \in \set{2, 4, 6, \cdots}$, the recursive process will stop at the nodes of leaves when $n -m = 2$ or $n-m = 0$.

\subsubsection{Space Complexity of BBTRA}

The space complexity of computing $\Radipoly{n}{m}(\rho)$  is determined by the memories needed. For $n = m$, we just need one variable to record the value of $\Radipoly{n}{m}(\rho)$ once the radius $\rho$ and integer $m$ is given, thus the corresponding space complexity is $\mathscr{O}(1)$. However, for $n\neq m$, all of the nodes in the binary tree determined by the recursive process will be kept for computing $\Radipoly{n}{m}(\rho)$. Hence the storage is determined by the number of nodes of the BBT. As demonstrated by \Fig \ref{fig-recursive}, the number of nodes is 
\begin{equation}
\scc{\ProcName{CalcRadiPolyBBTRA}}{n,m} = 2^k -1 = \BigO{2^{\frac{n-m}{2}}}.
\end{equation}
In consequence, the space complexity $\BigO{2^{n-m}}$ is exponential, which increases very fast with the difference  $n-m$. In the sense of space complexity or memory consumption, the recursive scheme is not recommended  for large value of $n-m$. Particularly, for $n = m$ and $n = m+2$, the space  complexity is $\BigO{1}$.

\subsubsection{Time Complexity of BBTRA}

The crucial result in our development of time complexity will be the following fundamental equation 
\begin{equation}
\tcc{\cpvar{op}}{\Radipoly{n-1}{m+1}} = \tcc{\cpvar{op}}{\Radipoly{n-2}{m}}, \quad \cpvar{op}\in \set{*, +} 
\end{equation}
for the non-leaf nodes of the BBT, which simplifies the analysis of computational complexity in a surprising way. Let
\begin{equation}
f_{\cpvar{op}}(n,m) = \tcc{\cpvar{op}}{\Radipoly{n}{m}},
\end{equation} 
then we have 
\begin{equation}
\left\{
\begin{array}{ll}
f_+(n-1, m+1) = f_+(n-2,m), & n -m \ge 4;\\
f_*(n-1, m+1) = f_*(n-2,m), & n -m \ge 4.
\end{array}
\right.
\end{equation}
According to \eqref{eq-tcc-leaf-type-A} and \eqref{eq-tcc-leaf-type-B}, for $m\ge 1$ we can obtain
\begin{equation}
\left\{
\begin{array}{ll}
f_+(n,m)  = 3, & m\neq 0, n-m = 2; \\
f_+(n, m) = 0, & m\neq 0, n-m = 0.
\end{array}
\right.
\end{equation}
and
\begin{equation}
\left\{
\begin{array}{ll}
f_*(n, m) = 3 + \mceil{\log_2 m}, & m\neq 0, n-m = 2; \\
f_*(n, m) = \mceil{\log_2 m}, & m\neq 0, n-m = 0.
\end{array}
\right.
\end{equation}
For $m=0$, simple algebraic calculations show that
\begin{equation}
\left\{
\begin{array}{ll}
f_+(n,m) = 3, &m = 0,  n-m = 2; \\
f_+(n,m) = 0, &m = 0,  n-m = 0.
\end{array}
\right.
\end{equation}
and 
\begin{equation}
\left\{
\begin{array}{ll}
f_*(n, m) = 3, &m = 0,  n-m = 2; \\
f_*(n, m) = 0, &m = 0,  n-m = 0.
\end{array}
\right.
\end{equation}
For the addition operations, equation \eqref{eq-R-BBT} implies that
\begin{equation}
\begin{split}
\tcc{+}{\Radipoly{n}{m}} 
&= \tcc{+}{\rho F_1\Radipoly{n-1}{m+1}} + 
\tcc{+}{F_2\Radipoly{n-2}{m}} + 1\\
&=  \tcc{+}{\Radipoly{n-1}{m+1}} + \tcc{+}{\rho \frac{2n}{n-m}}
 + \tcc{+}{1-F_1} +  \tcc{+}{\Radipoly{n-2}{m}} + 1 \\
&= 2\tcc{+}{\Radipoly{n-1}{m+1}} + 3
\end{split}
\end{equation}
Thus
\begin{equation}
f_+(n,m) = 2f_+(n-1,m+1)+ 3, \quad n-m\ge 4
\end{equation}
Let 
\begin{equation}
G^+_p = f_+(n,m) = f_+(n-m), \quad p=n-m\in \set{0, 2, 4, \cdots}
\end{equation}
then 
\begin{equation}
G^+_{p-2} = f_+(n-1,m+1) = f_+(n-2, m),
\end{equation}
hence we have the discrete difference equation
\begin{equation} \label{eq-Gp+}
G^+_p = 2 G^+_{p-2} + 3, \quad G^+_2 = 3, G^+_0 = 0.
\end{equation}
According to \eqref{eq-Gp-de} and \eqref{eq-ans-Gp-de}, the solution to \eqref{eq-Gp+} must be 
\begin{equation}
G^+_p 
= 2^{\frac{p}{2}-1}G^+_2 + 3(2^{\frac{p}{2}-1}-1)
= 3(2^{\frac{p}{2}}-1), \quad p = 0, 2, 4, \cdots
\end{equation}
Hence we have time flops for the addition operation
\begin{equation}
f_+(n,m) = \tcc{+}{\Radipoly{n}{m}} = 3(2^{\frac{n-m}{2}}-1), \quad n-m= 0, 2, 4, \cdots
\end{equation}
Similarly, for the multiplication operations, we have
\begin{equation}
\begin{split}
\tcc{*}{\Radipoly{n}{m}} 
= \tcc{*}{\Radipoly{n-1}{m+1}} + 
\tcc{*}{\Radipoly{n-2}{m}} 
 + \tcc{*}{\rho F_1} + \tcc{*}{F_2} + 2  
= 2\tcc{*}{\Radipoly{n-1}{m+1}} + 5.
\end{split}
\end{equation}
In this way, we immediately obtain the following discrete difference equation
\begin{equation} \label{eq-Gp*}
G^*_p = 2G^*_{p-2} + 5, \quad p = 2, 4, 6, \cdots
\end{equation}
and the initial value
\begin{equation} \label{eq-G2*}
G^*_2 = 
\left\{
\begin{array}{ll}
\mceil{\log_2 m}, & m \ge 1; \\
3, & m = 0.
\end{array}
\right.
\end{equation}

With the help of \eqref{eq-Gp-de}, \eqref{eq-ans-Gp-de}, \eqref{eq-Gp*} and \eqref{eq-G2*}, we can deduce that 
\begin{equation}
\begin{split}
G^*_p 
= 2^{\frac{p}{2}-1}G^*_2 + 5\left(2^{\frac{p}{2}-1}-1\right) 
=\left\{
\begin{array}{ll}
8\cdot 2^{\frac{n-m}{2}-1} - 5, & m = 0; \\
\left(\mceil{\log_2 m} + 8\right)2^{\frac{n-m}{2}-1} - 5, & m \ge 1.
\end{array}
\right.
\end{split}
\end{equation} 
Therefore,
\begin{equation}
\begin{split}
f_*(n,m) 
= \tcc{*}{\Radipoly{n}{m}} 
=\left\{
\begin{array}{ll}
\left(\mceil{\log_2 m} + 8\right)2^{\frac{n-m}{2}-1} - 5, & m\neq 0, n - m \ge 2;\\
\mceil{\log_2 m}, & m\neq 0, n-m = 0;\\
1, & m = 0, n-m = 0 \\
8 \cdot 2^{\frac{n-m}{2}-1} - 5, & m = 0, n - m \ge 2\\
\end{array}
\right.
\end{split}
\end{equation}
The total time flops for computing $\Radipoly{n}{m}(\rho)$ is
\begin{equation}
\begin{split}
\Tf(\Radipoly{n}{m}) 
= \tcc{*}{\Radipoly{n}{m}} +  \tcc{+}{\Radipoly{n}{m}} 
=\left\{
\begin{array}{ll}
\left[7 + \frac{\mceil{\log_2 m}}{2}\right]2^{\frac{n-m}{2}} - 5, & m\neq 0, n-m \ge 2 \\
\mceil{\log_2 m}, & m\neq 0, n-m = 0 \\
7 \cdot 2^{\frac{n-m}{2}} - 5, & m = 0, n-m \ge 2 \\
1, & m = 0, n-m = 0 \
\end{array}
\right.
\end{split}
\end{equation} 
which can be expressed with the big-O notation as follows
\begin{equation} \label{eq-tcc-Rnm}
\begin{aligned}
\Tf(\Radipoly{n}{m})
=\left\{
\begin{array}{ll}
\BigO{\log_2 \abs{m}2^{\frac{n-\abs{m}}{2}}}, & m\neq 0, n- \abs{m}\ge 4;\\ 
\BigO{\log_2 \abs{m}}, & m\neq 0,  n - \abs{m} =0, 2; \\
\BigO{2^{\frac{n-\abs{m}}{2}}}, & m = 0, n- \abs{m}\ge 0.\\ 
\end{array}
\right.
\end{aligned}
\end{equation}
Clearly, in the sense of time complexity or time consumption, the recursive scheme is not recommended  for the large value of $n-\abs{m}$.

\section{Iterative Scheme for Computing the Zernike Radial Polynomials}
\label{sec-iterative}

For the recursive \Algr \ref{alg-Rnm-recur}, its space complexity is exponential and the time complexity varies with the configuration of double indices $n$ and $m$ according to \eqref{eq-tcc-Rnm}. For large value of $n-\abs{m}$, it is necessary for us to reduce the computational complexity with iterative scheme instead of recursive scheme for computing $\Radipoly{n}{m}(\rho)$.  

\subsection{Computing the Radial Polynomials Iteratively}

The key issue of iterative computation lies in two points: firstly, we set the initial condition in an iterative scheme with the stopping condition in the corresponding recursive scheme; then we iterate by updating the state of nodes of interest via the primitive  recursive formulae in a reverse direction. For the recursive problem with a single integer parameter, say $n$, this can be done easily since there are just two directions for the variation of $n$ for the 1-dim problems. However, for the recursive problem with double integer parameters, say $n$ and $m$, the iteration may be trouble because there is no simple direction for the  updating process. Fortunately, for computing the ZRP, it is possible to iterate simply with the structure of BBT. Actually, what we should do is just accessing the nodes in BBT from the leaf nodes to the root node.

\begin{figure*}[h]
\centering
\includegraphics[width=1.0\textwidth]{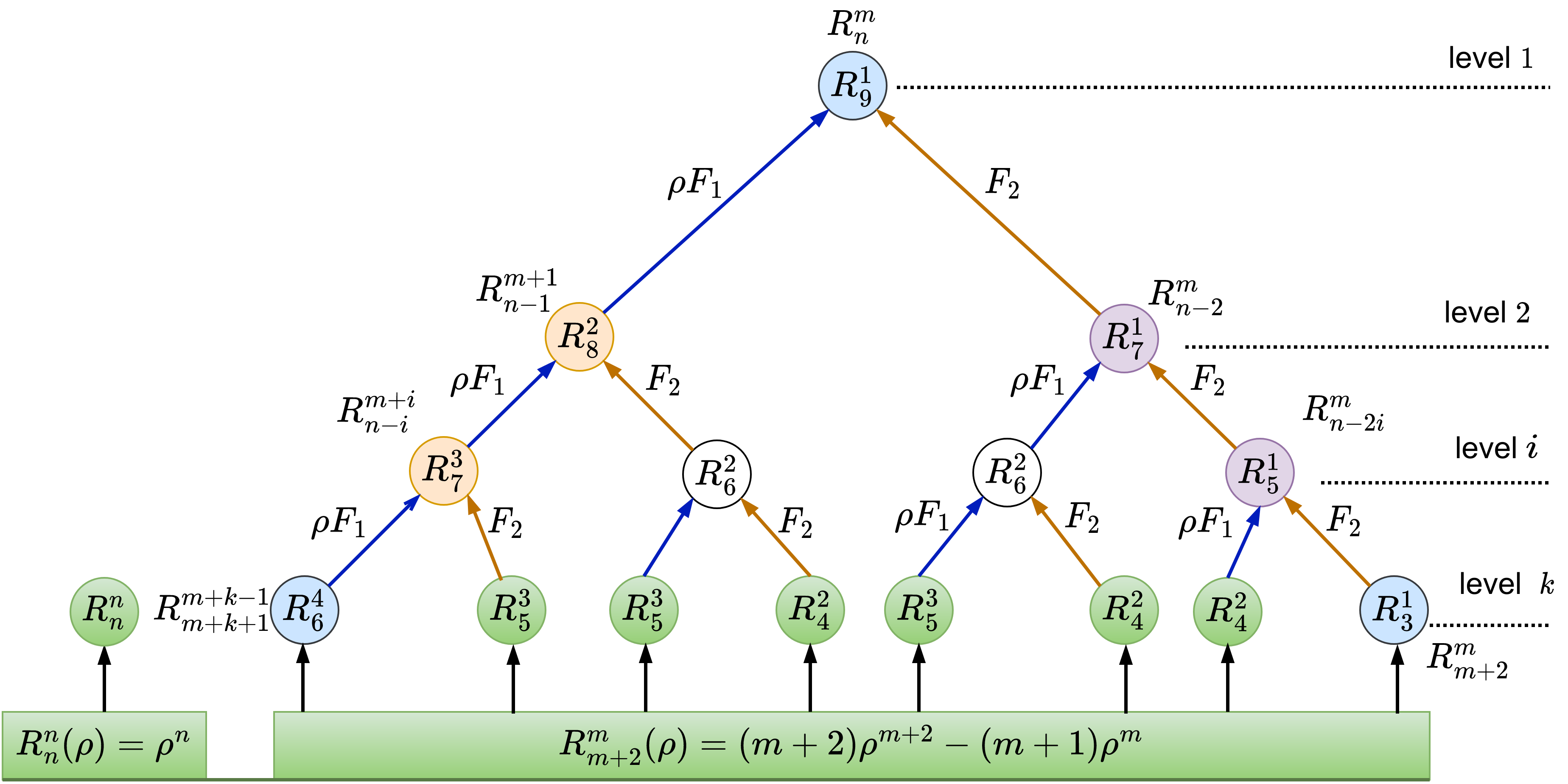} 
\caption{Principle of Iterative Scheme for Computing $\Radipoly{n}{m}(\rho)$ when $n\neq \abs{m}$}
\label{fig-BBT-principle-iter}
\end{figure*}

\Fig \ref{fig-BBT-principle-iter} shows the iterative process intuitively. For $m\ge 0$ and $n-m\ge 2$, we start with the initial condition specified by 
\eqref{eq-R-sc} in the leaves in the $k$-th level, then compute the value of the upper nodes  in the $(k-1)$-th level according to \eqref{eq-Rnm-recu}, finally update the state of the nodes. The formula in \eqref{eq-R-BBT} can be used for both recursion and iteration, which depends on the concrete 
usages. Moreover, it is obvious that there are some repetitive nodes in \Fig \ref{fig-BBT-principle-iter}, which implies that we can save their calculation to accelerate the computation and reduce memories needed. 

\begin{figure*}[h]
\centering
\includegraphics[width=1.0\textwidth]{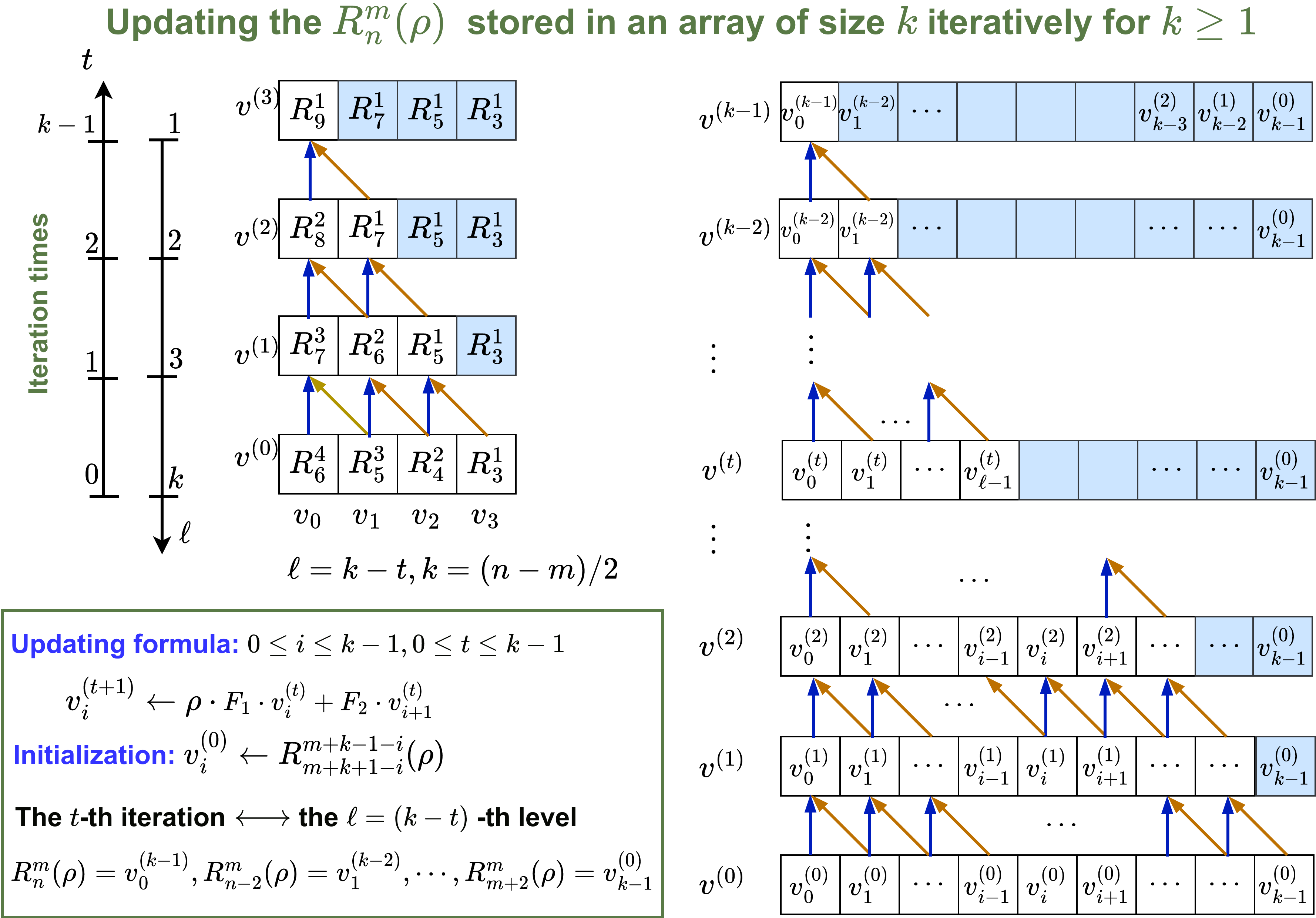} 
\caption{Illustration of Iterative Algorithm for Computing $\Radipoly{n}{m}(\rho)$}
\label{fig-BBT-array-iter}
\end{figure*}

\Fig \ref{fig-BBT-array-iter} demonstrates the simplification of the BBT with a 1-dim array to store the nodes of interest. By removing the redundant  nodes,  there are just $k$ independent leaves in the $k$-th level of the BBT. We can use a sequence $\cpvar{v}=\seq{v_0, v_1, \cdots, v_{k-1}}$ to store the initial values with the expression
\begin{equation}
\begin{split}
v_i 
= \Radipoly{m+k+1-i}{m+k-1-i}(\rho) 
=\rho^{m+k-i-1}[(m+k-i+1)\rho^2-(m+k-i)]
\end{split}
\end{equation} 
for $0\le i \le k-1$.
With the help of \eqref{eq-Rnm-recu}, we can deduce that the iterative formulae for updating must be
\begin{equation}
v_i^{(t+1)} = \rho F_1 v_i^{(t)} + F_2 v_{i+1}^{(t)}, \quad i, t\in \set{0, 1, \cdots, k-1}
\end{equation}
where $t$ denotes the $t$-th iteration. It should be noted that the direction of iteration is from the bottom to top such that 
\begin{equation}
t = k - \ell.
\end{equation}
In other words, the variable $t$ for iteration can be replaced by the index $\ell$ for the levels of the BBT. When $t = k-1$, the iterative process stops and we have
\begin{equation}
\Radipoly{n}{m}(\rho) = v_0^{(k-1)}.
\end{equation}
When the iteration ends, we obtain the byproduct stored in the sequence $\cpvar{v}$:
\begin{equation}
\begin{split}
\seq{v_0^{(k-1)}, v_1^{(k-2)}, \cdots, v_{k-1}^{(0)}}
= \seq{\Radipoly{n}{m}(\rho), \Radipoly{n-2}{m}(\rho), \cdots, \Radipoly{m+2}{m}(\rho)},
\end{split}
\end{equation}
in which only $v_0^{(k-1)}$ is the computation result of $\Radipoly{n}{m}(\rho)$ and the others are discarded. 

In the left part of \Fig \ref{fig-BBT-array-iter}, the computation of $\Radipoly{9}{1}(\rho)$ is illustrated intuitively where $n = 9, m = 1$ and $ k = 4$. As a generalization, the right part of \Fig \ref{fig-BBT-array-iter} demonstrates the iterative process of computing $\Radipoly{n}{m}(\rho)$ for $n-m\ge 2$ and $m\ge 0$. 

\subsection{Balanced Binary Tree Iterative Algorithm (BBTIA)}

Although the principle of converting the recursive process to its iterative counterpart is intuitive, it is still necessary to present the iterative algorithm clearly for the purpose of implementation with computer programs. \Algr \ref{alg-Rnm-iter} is the \textit{balanced binary tree iterative algorithm} (BBTIA)  for computing $\Radipoly{n}{m}(\rho)$ in an iterative way. We remark that the dynamic memory allocation in \Algr \ref{alg-Rnm-iter}
 is necessary because the parameter $k$ is determined by the input double indices $\mpair{n}{m}$. 

\begin{breakablealgorithm}
\caption{Compute the radial function $\Radipoly{n}{m}(\rho)$  iteratively} \label{alg-Rnm-iter}
\begin{algorithmic}[1]
\Require Born-wolf double indices $\mpair{n}{m}$, variable $\rho\in[0, 1]$
\Ensure the value of $\Radipoly{n}{m}(\rho)$
\Function{CalcRadiPolyBBTIA}{$n$, $m$, $\rho$}
\State $m \gets \abs{m}$;  \quad // $\Radipoly{n}{m}(\rho) = \Radipoly{n}{-m}(\rho)$.
\State $p\gets n-m$;
\If{$p=0$}
\State \Return $\ProcName{CalcLeafNodeTypeA}(\rho,m)$;
\EndIf
\If{$p=2$}
\State \Return $\ProcName{CalcLeafNodeTypeB}(\rho, m)$;
\EndIf
\State $ k \gets p/2$; \quad // $k = (n-m)/2$
\State Allocate memories for the array $\cpvar{v}=\seq{v_0, v_1, \cdots, v_{k-1}}$;
\For{$i\in\seq{0,1,\cdots,k-1}$}
\State $v_i = \ProcName{CalcLeafNodeTypeB}(\rho, m+k-1-i)$;
\EndFor
\For{$\ell \in \seq{k-1, k-2, \cdots, 1}$}
	\State $\scrd{n}{start} \gets  n - (\ell -1)$; \quad // number $n$ for the most left node in level $\ell$.
	\State $\scrd{m}{start} \gets  m + (\ell -1)$; \quad // number $m$ for the most left node in level $\ell$.
	\For{$i\in \seq{0,1,\cdots,\ell-1}$} 	\quad // compute the node $R_{\scrd{n}{pos}}^{\scrd{m}{pos}}(\rho)$ in the level $\ell$
    	\State $\scrd{m}{pos} \gets \scrd{m}{start} - i, \scrd{n}{pos} \gets \scrd{n}{start} - i$;
    	\State $F_1 \gets 2\scrd{n}{pos}/(\scrd{n}{pos} -\scrd{m}{pos}), F_2 \gets 1- F_1$;
    	\State $v_i \gets \rho F_1 v_i + F_2 v_{i+1}$;  \quad// Updating the $i$-th value $v_i$ in the $(k -\ell)$-th level;
    \EndFor
\EndFor
\State $R\gets v_0$;
\State Release the memories for the sequence \cpvar{v};
\State \Return $R$;
\EndFunction
\end{algorithmic}
\end{breakablealgorithm}

\subsection{Computational Complexity of BBTIA}

\subsubsection{Space Complexity of BBTIA}

For $n=m$ and $n=m+2$, equation \eqref{eq-R-sc}  can be used and only one variable is needed to store the value of  $\Radipoly{n}{m}(\rho)$, hence the space complexity is $\mathscr{O}(1)$.
For $n\ge m + 4$, the sequence $\cpvar{v} = \seq{v_0, v_1, \cdots, v_{k-1}}$ is used to realize the computation of $\Radipoly{n}{m}(\rho)$. The number of auxiliary variables (i.e., $p, k, \scrd{n}{start}, \scrd{m}{start}, F_1, F_2$ and $R$) is constant, hence memory assumption is dominated by the length $k = n-\abs{m}$ of the sequence \cpvar{v}. In summary,  the space complexity for computing $\Radipoly{n}{m}(\rho)$ is 
\begin{equation}
\scc{\ProcName{CalcRadiPolyBBTIA}}{n,m} = \BigO{n-\abs{m}}, \quad \abs{m}\le n. 
\end{equation}

\subsubsection{Time Complexity of of BBTIA}

\begin{table*}[htb]
\centering
\caption{Counting the time flops in \Algr \ref{alg-Rnm-iter} for $m\ge 1$.}
\begin{tabular}{ccccl}
\hline
\cpvar{Line}   & $\tcc{+}{\cpvar{Line}}$ &  $\tcc{*}{\cpvar{Line}}$ & Loop Counting  &  Remark   \\
\hline
5 &  0   &  $\mceil{\log_2 m}$   & 1 & for leaf node $\Radipoly{m}{m}$ \\  
8 & $3$   &  $3+\mceil{\log_2 m}$ &  1 & for leaf node $\Radipoly{m+2}{m}$    \\
10 &   0  & 1               & 1        &   for root node $\Radipoly{n}{m}$ \\
13 &   $3$  & $3+\mceil{\log_2(m+k-1-i)}$  & $\displaystyle \sum^{k-1}_{i=0}$  & for root node $\Radipoly{n}{m}$ \\
16 &    2&    0                &  $\displaystyle \sum^{k-1}_{\ell =1} $ & for root node $\Radipoly{n}{m}$  \\
17 &    2&       0             & $\displaystyle \sum^{k-1}_{\ell =1} $  &  for root node $\Radipoly{n}{m}$  \\
19 &    2&     0               & $\displaystyle \sum^{k-1}_{\ell =1}\sum^{\ell-1}_{i=0} $ & for root node $\Radipoly{n}{m}$ \\
20 &  2   & 2                 & $\displaystyle \sum^{k-1}_{\ell =1}\sum^{\ell-1}_{i=0} $ & for root node $\Radipoly{n}{m}$ \\
21 &  1   &  3                & $\displaystyle \sum^{k-1}_{\ell =1}\sum^{\ell-1}_{i=0} $ &  for root node $\Radipoly{n}{m}$ \\
\hline 
\end{tabular}
\end{table*}

For $n - m = 0$ and $n - m = 2$, the time computational complexity is determined by the Line 5 and Line 8 in \Algr \ref{alg-Rnm-iter} respectively. The time complexity for the leaf nodes of type $\ProcName{X}$ is 
\begin{equation}
\tcc{\ProcName{CalcLeafNodeTypeX}}{n, m} = 
\left\{
\begin{array}{ll}
\BigO{1}, & n - m = 0; \\
\BigO{\log_2 m}, & m\neq 0, n - m = 2;
\end{array}
\right.
\end{equation}
according to \eqref{eq-tcc-power} where X is A or B.

For $m> 0$ and $n-m \ge 4$, the time computational complexity is determined by the operations in the loops involved in Lines $1\sim 27$ in \Algr \ref{alg-Rnm-iter}.  Let $I = \set{10, 13, 16, 17, 19, 20, 21}$ and ignore other line which the  time cost is not important for complexity analysis, we can obtain the addition flops  
\begin{equation}
\begin{split}
\tcc{+}{\Radipoly{n}{m}}
= \sum_{i \in I}\tcc{+}{\cpvar{State}[i]}
= \frac{5}{2}k^2 + \frac{9}{2}k - 4
= \frac{5}{8}(n-m)^2 + \frac{9}{4}(n-m)  -4
\end{split}
\end{equation}  
since $k = (n-m)/2$ and multiplication flops 
\begin{equation}
\tcc{*}{\Radipoly{n}{m}}
= \sum_{i \in I}\tcc{*}{\cpvar{State}[i]}
= \frac{5}{2}k^2 + \frac{1}{2} k + 1 
+ \sum^{k-1}_{i=0}\mceil{\log_2(m+k-1-i)}
\end{equation}  
respectively. 
Consequently, we have 
\begin{equation}
\begin{split}
\tcc{\ProcName{CalcRadiPolyBBTIA}}{n,m} 
= \BigO{k^2} 
= \BigO{(n-m)^2}, \quad m> 0, n-m\ge 4.
\end{split} 
\end{equation}
As an illustration, for $\Zernpoly{94}{}(\rho,\theta)= \sqrt{28} \Radipoly{13}{1}(\rho)\sin \theta$ we have $n = 13, m = 1, k = (n-m)/2=6$, thus we can obtain
\begin{align*}
\tcc{*}{\Radipoly{13}{1}} 
= \frac{5}{2}\cdot 5^2 + \frac{5}{2} + 1+ \sum^4_{i=0}\mceil{\log_2(6-i)} = 74
\end{align*} 
Particularly, for $m = 0$, we can obtain \Tab \ref{tab-flop-Rn0}. The time flops are also determined by the operations in the loops involved on Line $1\sim 27$ in \Algr \ref{alg-Rnm-iter}. With the help of \Tab \ref{tab-flop-Rn0}, an argument similar to the process for finding $\tcc{\cpvar{op}}{\Radipoly{n}{m}}$ where $m>0$ shows that 
\begin{equation}
\tcc{\cpvar{op}}{\Radipoly{n}{0}} 
= \BigO{(n-m)^2}, \quad \cpvar{op}\in \set{*, +}.
\end{equation} 
since $k=(n-m)/2$ for non-negative $m$.
\begin{table*}[htb]
\centering
\caption{Counting the time flops in \Algr \ref{alg-Rnm-iter} for $m= 0$.} 
\label{tab-flop-Rn0}
\begin{tabular}{ccccl}
\hline
\cpvar{Line}  & $\tcc{+}{\cpvar{Line}}$ &  $\tcc{*}{\cpvar{Line}}$ & Loop Counting  &  Remark   \\
\hline
5 &  0   &  0   & 1 & for the leaf node $\Radipoly{0}{0}(\rho)$ \\  
8 & $3$   &  $3$ &  1 & for the leaf node $\Radipoly{2}{0}$    \\
10 &   0  & 1               & 1        &   for the root node $\Radipoly{n}{0}(\rho)$ \\
13 &   $3$  & $2$  & 1  & for the leaf $\Radipoly{2}{0}(\rho)$ \\
   &   $3$  & $2 + \mceil{\log_2 (k-1-i)}$  & $\displaystyle \sum^{k-2}_{i=0}$  & for the leaf $\Radipoly{k+1-i}{k-1-i}(\rho)$ \\
16 &    2&    0                &  $\displaystyle \sum^{k-1}_{\ell =1} $ & for the root node $\Radipoly{n}{0}(\rho)$  \\
17 &    2&       0             & $\displaystyle \sum^{k-1}_{\ell =1} $  &  for the root node $\Radipoly{n}{0}(\rho)$  \\
19 &    2&     0               & $\displaystyle \sum^{k-1}_{\ell =1}\sum^{\ell-1}_{i=0} $ & for the root node $\Radipoly{n}{0}(\rho)$ \\
20 &  2   & 2                 & $\displaystyle \sum^{k-1}_{\ell =1}\sum^{\ell-1}_{i=0} $ & for the root node $\Radipoly{n}{0}(\rho)$ \\
21 &  1   &  3                & $\displaystyle \sum^{k-1}_{\ell =1}\sum^{\ell-1}_{i=0} $ &  for the root node $\Radipoly{n}{0}(\rho)$ \\
\hline 
\end{tabular}
\end{table*}

As a comparison, \Tab \ref{tab-complexity-Rnm} shows the computational complexity with the $\mathscr{O}(\cdot)$ notation for computing the ZRP $\Radipoly{n}{m}(\rho)$.

\begin{table*}[htb] 
\centering
\caption{Computational Complexity of Computing $\Radipoly{n}{m}(\rho)$}
\label{tab-complexity-Rnm}
\begin{tabular}{llll}
\hline
Type of Algorithm  & Double Indices $\mpair{n}{m}$  & Space Complexity & Time Complexity \\
\hline
Recursive algorithm & $n-\abs{m}\in \set{0, 2}, m\neq 0$  & $\BigO{1}$ &  $\BigO{\log_2 \abs{m}}$ \\
\ProcName{CalcRadiPolyBBTRA}  & $n-\abs{m}\in \set{0, 2}, m = 0$  & $\BigO{1}$ &  $\BigO{1}$ \\
   & $n- \abs{m}\ge 4, \abs{m} =1$ & $\BigO{2^{\frac{n-\abs{m}}{2}}}$ & $\BigO{2^{\frac{n-\abs{m}}{2}}}$  \\
                    & $n- \abs{m}\ge 4, \abs{m}> 1$&  $\BigO{2^{\frac{n-\abs{m}}{2}}}$ & $\BigO{(\log_2 \abs{m}) 2^{\frac{n-\abs{m}}{2}}}$  \\   
\hline
Iterative algorithm & $n -\abs{m} \in \set{0,2}, m \neq 0$  & $\BigO{1}$ &  $\BigO{\log_2 \abs{m}}$ \\
\ProcName{CalcRadiPolyBBTIA} & $n -\abs{m} \in \set{0,2}, m = 0$  & $\BigO{1}$ &  $\BigO{1}$ \\
    & $n -\abs{m} \ge 4$        & $\BigO{n-\abs{m}}$ & $\BigO{(n-\abs{m})^2}$ \\ 
\hline
\end{tabular}
\end{table*}

\section{Verification and Validation} \label{sec-v-and-v}

We implemented the novel recursive and iterative algorithms with the C programming language and compared the running time with Prata-Rusch's and Shakibaei-Paramesran's recursive methods. As an illustration, we have tested the practical running time for $n = 28$ and $n= 29$ respectively (where $m \in \set{n, n-2, \cdots, n-2\mfloor{n/2}}$) via an average value by repeating the algorithms of interest for 10 times.

\begin{figure*}[htb]
\centering
\subfigure[Running time for $\Radipoly{28}{m}(\rho)$]{
\includegraphics[width=0.45\textwidth]{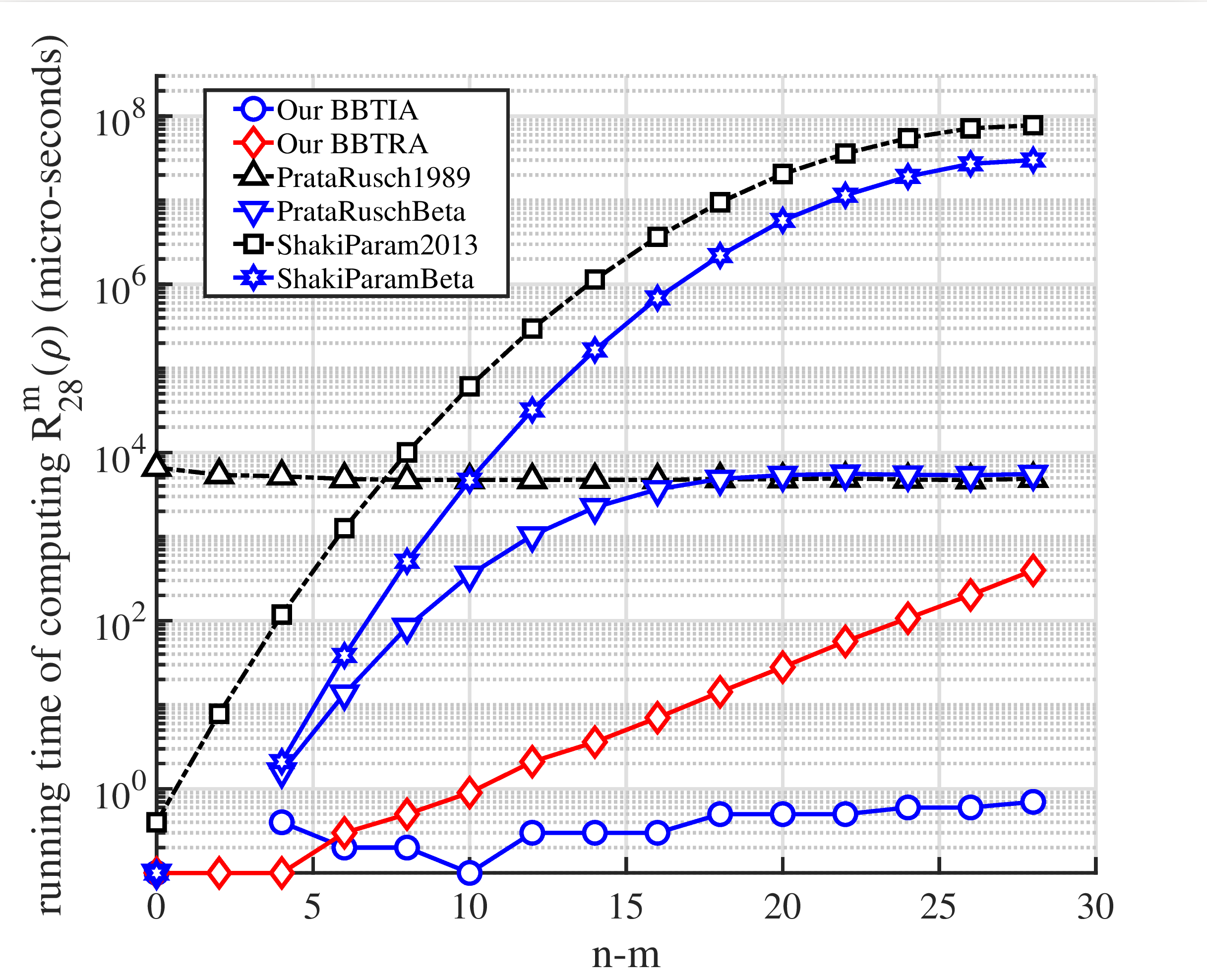} 
}
\subfigure[Running time for $\Radipoly{29}{m}(\rho)$]{
\includegraphics[width=0.45\textwidth]{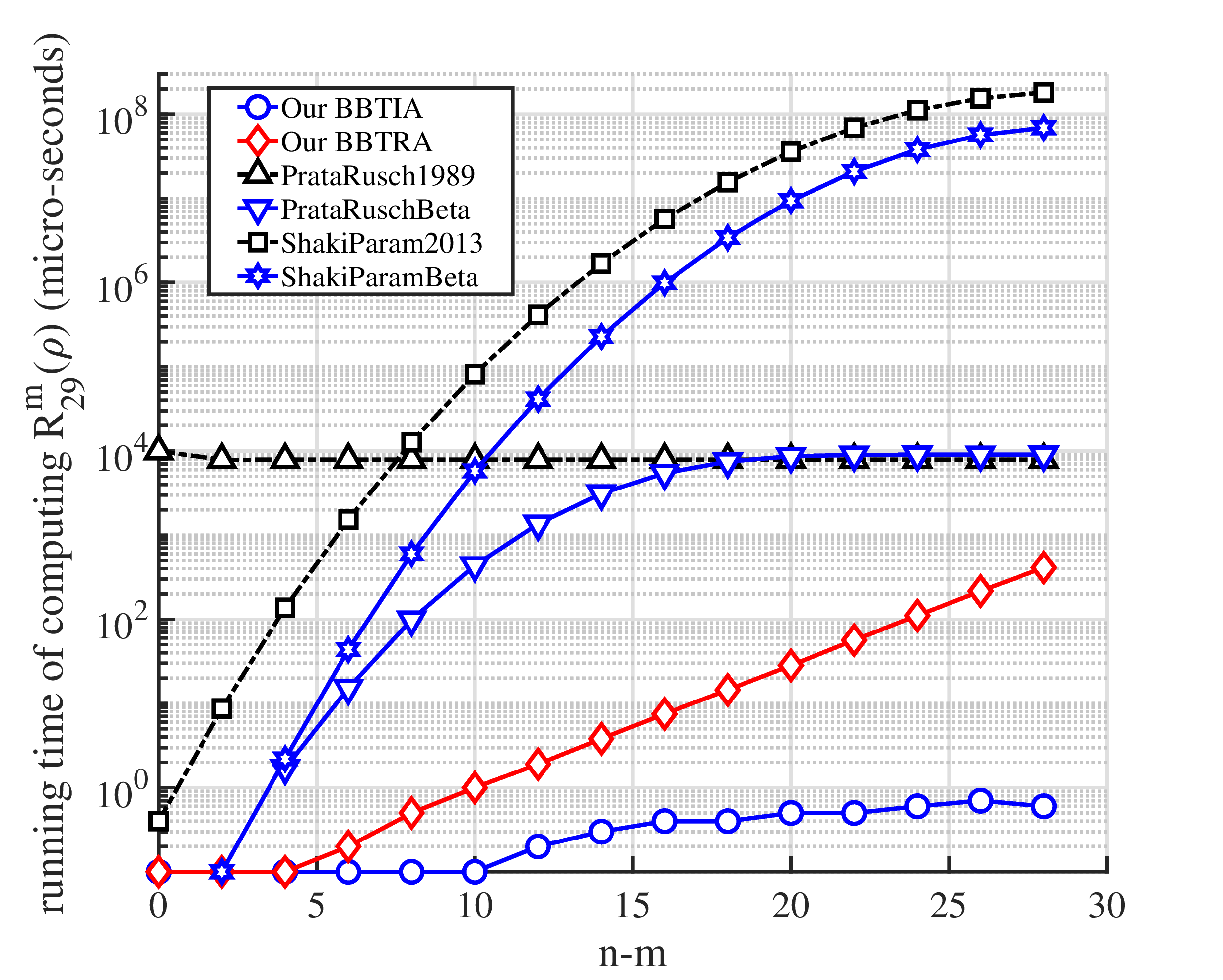} 
}
\caption{Diagrams of running time of computing $\Radipoly{n}{m}(\rho)$ the difference of $n$ and $m$ for $n\in \set{28, 29}$ and $m\ge 0$.}
\label{fig-runtime-comp}
\end{figure*} 

As shown in \Fig \ref{fig-runtime-comp} from top to bottom, there are six curves of running time varying with the difference $n-m$. Here we give some necessary interpretations:
\begin{itemize}
\item the blue curve marked by circle $\circ$ and legend "\texttt{Our BBTIA}" corresponds to the running time of our BBTIA (\Algr \ref{alg-Rnm-iter}), which is at the bottom in the subfigures;
\item the red curve marked by diamond and legend "\texttt{Our BBTRA}" corresponds to the running time of our BBTRA (\Algr \ref{alg-Rnm-recur}), which looks like a straight line in the semilog diagram starting from $n-\abs{m}=4$ in the subfigures since for $n-\abs{m}\ge 4$ the time complexity is exponential;
\item the dashed curve marked by upward-facing triangle and legend "\texttt{PrataRusch1989}" corresponds to the running time of the primitive Prata-Rusch
recursive scheme with the stopping condition \eqref{eq-R-sc-PR-v1}, which also looks like a horizontal line in the semilog diagram; it is over the straight line for our recursive \Algr \ref{alg-Rnm-recur}, which implies the time computational complexity is much higher than the exponential complexity;
\item the blue curve marked by downward-facing triangle and legend "\texttt{PrataRuschBeta}" corresponds to the running time of the improved Prata-Rusch
recursive scheme with the stopping condition \eqref{eq-R-sc},
which reduces the time consumption for $n-\abs{m}\le 18$ and remains the same with the counterpart of the  primitive Prata-Rusch recursive scheme;
\item the dashed curve marked by square and legend "\texttt{ShakiParam2013}" corresponds to the running time of the primitive Shakebaei-Paramesran
recursive scheme with the stopping condition \eqref{eq-R-sc-SP}
\item the blue curve marked by hexagram and legend "\texttt{ShakiParamBeta}" corresponds to the running time of the improved Shakebaei-Paramesran
recursive scheme with the stopping condition \eqref{eq-R-sc}.
\end{itemize} 
\Fig \ref{fig-runtime-comp}  shows clearly that 
\begin{itemize}
\item the running time for our iterative algorithm is the shortest since its time computational complexity is just $\BigO{(n-\abs{m})^2}$ and
\item  the semilog curve for the running time for our recursive algorithm looks like a straight line, which 
corresponds to the exponential complexity and coincides with the theoretical analysis above very well.
\end{itemize}
\Fig \ref{fig-runtime-comp} also illustrates a quantitative result that for $n -m = 28$ (i.e., $n=28, m= 0$ or $n = 29, m = 1$): the worst case of time consumption for computing $\Radipoly{n}{m}(\rho)$ is about $10^7\sim 10^8$ microseconds ($10\sim 100$ seconds) with recursive  schemes, our BBTIA works very well and the time consumption is about $1$ microseconds, which implies that it is suitable for real-time applications.  

The trends of curves in \Fig \ref{fig-runtime-comp} show that the time consumed increases with the growth of $n-m$ when $n-m \ge 4$ and remains constant when $n-m \in \set{0, 2}$. In \cite{Shakibaei-2013}, the conclusion that the time complexity  depends only on $n$ is arguable since the smaller the $n-m$, the faster the algorithms are.   
The computational complexity of recursive algorithms by Shakibaei-Paramesran and Prata-Rusch is much higher than the exponential complexity since the curve of time consumption is over the line specified by our recursive algorithm with exponential complexity. For our iterative algorithm with square complexity $\BigO{(n-m)^2}$, the running time is within 1 micro-second, which implies that it is suitable for real-time applications. 

The data for \Fig \ref{fig-runtime-comp} can be downloaded from the GitHub website:
 \blue{\url{https://github.com/GrAbsRD/Zernike/tree/RadialPolynomialRunTime}}.
It should be pointed out that the configuration of our testing platform is as follows: 
Debian GNU/Linux 11 (bullseye) OS; Memory, 64GB RAM; Processor, AMD$^\circledR$ Ryzen 7 5800$\times$8-core processor$\times$16; Compiler, gcc-10.2.1 20210110 (Debian 10.2.1-6). For other computational platform, the running time may be different but the trends of the curves should be similar and the straight line for our recursive \Algr \ref{alg-Rnm-recur} will still remain a straight line    in the semilog diagram which illustrates the exponential complexity $\BigO{2^{\frac{n-\abs{m}}{2}}}$ clearly.
Furthermore, the time consumption of our iterative \Algr \ref{alg-Rnm-iter} will stay at the bottom since its $\BigO{(n-\abs{m})^2}$ computational complexity is the lowest.

\section{Conclusion} \label{sec-summary}

The numerical computation of ZRP is a challenging problem due to the stability and the computational complexity. Our formulae  
\eqref{eq-Rnm-recu} for computing $\Radipoly{n}{m}(\rho)$ owns the following advantages: 
\begin{itemize}
\item firstly it is capable of computing the value of ZRP with high precision by avoiding computing the  high order power functions $\rho^\ell$ where $\ell$ is a large positive integer;  
\item secondly it has a simple algebraic structure for understanding, remembering  and applications;
\item thirdly it is symmetric which leads to a balanced binary tree for recursion and convenient theoretic analysis of computational complexity;
\item fourthly the conversion of recursive process to the  iterative version is easy and the quadratic complexity refreshes the state-of-the-art of the computational complexity for computing $\Radipoly{n}{m}(\rho)$; and
\item finally it leads to a stable computation process since there is no singularity in the expression.
\end{itemize}

The BBTRA proposed is the fastest recursive algorithm for computing the ZRP when compared with the other available recursive algorithms. The BBTIA is the first iterative algorithm  for computing the ZRP which has the quadratic time complexity and it is suitable for real-time applications. The high precision can be achieved automatically due to the recursive and iterative property for computing polynomials. In the sense of STEM education, the connection of the BBT and ZRP exhibits the beauty and applications of discrete mathematical structure behind the engineering problem, which is worthy of introducing to the college students, computer programmers and optics engineers.

\section*{Code and Data Availability Statement}

The code and the data for \Fig \ref{fig-runtime-comp} can be downloaded from the GitHub websites 
 \blue{\url{https://github.com/GrAbsRD/Zernike/tree/RadialPolynomialRunTime}} 
and 
\blue{\url{https://github.com/GrAbsRD/Zernike}} 
respectively.

\section*{Acknowledgment}

This work was supported in part by the Hainan Provincial
Natural Science Foundation of China under Grant 
2019RC199 and in part by the National Natural Science Foundation of China under Grant 62167003.

\appendix

\section{Notations for Integers and Computational Complexity  }

\subsection{Integers and Lower Integer for Real Number}
For an integer $n\in \mathbb{Z} = \set{0, \pm 1, \pm 2, \cdots}$, it is even if and only if $2\mid n$, i.e., $2$ divides $n$ or equivalently  $n\equiv 0~(\mod 2)$; otherwise, it is odd if and only if $ 
2 \nmid n$ or equivalently $n \equiv 1~(\mod 2)$. 
The set of non-negative integers are denoted by
$\mathbb{Z}^+ = \set{0, 1, 2, \cdots}$. 
For any real number $x\in \mathbb{R}$, the maximal lower bound $n\in \mathbb{Z}$ such that $n \le x $ is called the floor of $x$ and it is denoted by
\begin{equation}
\mfloor{x} = n=\arg \max_{m} \set{m\in \mathbb{Z}: x\ge m}.
\end{equation}
Similarly,  the minimal upper bound $n\in \mathbb{Z}$ such that $ x \le n$ is called the ceiling of $x$ and it is denoted by
\begin{equation}
\mceil{x} = n=\arg \min_{m} \set{m\in \mathbb{Z}: x\le m}.
\end{equation}

\subsection{Time and Space Complexities} 
\label{app-subsec-cc}

Let
$\tcc{*}{\cpvar{expr}}$ and
$\tcc{+}{\cpvar{expr}} $
be the counting or times of multiplication and addition in some operation expression $\cpvar{expr}$. The \textit{time  complexity vector of computation} (TCVC)  for $\cpvar{expr}$ is defined by
\begin{equation}
\timcomp(\cpvar{expr}) = [\tcc{*}{\cpvar{expr}}, \tcc{+}{\cpvar{expr}}].
\end{equation}
Note that we just list two components of $\timcomp(\cpvar{expr})$ here since the subtraction and division can be treated as addition  and multiplication respectively for real numbers.
The \textit{time flops} (TF) \textit{flops}  \cite{Golub2013}, for computing $\cpvar{expr}$ is
\begin{equation}
\Tf(\cpvar{expr}) 
= \tcc{*}{\cpvar{expr}} + \tcc{+}{\cpvar{expr}}
\end{equation}
in which $\tcc{*}{\cpvar{expr}} $ and $\tcc{+}{\cpvar{expr}}$ are the time flops of multiplication and addition respectively.  

Similarly, we use $\Tc(\cpvar{expr})$ to denote the computation time for $\cpvar{expr}$. We also use $\Tc_*(\cpvar{expr})$ and $\Tc_+(\cpvar{expr})$ to represent the computation time for the multiplications and additions involved in $\cpvar{expr}$.
Given the time units $\tau_*$ and $\tau_+$ for multiplication and addition respectively, then  
\begin{equation}
\vec{\tau} = [\tau_*, \tau_+]
\end{equation}
will be the vector of time units. The time for computing $\cpvar{expr}$ can be represented by
\begin{equation}
\begin{split}
\Tc(\cpvar{expr})=& \braket{\timcomp(\cpvar{expr})}{\vec{\tau}} 
= \tcc{*}{\cpvar{expr}} \tau_* + \tcc{+}{\cpvar{expr}}\tau_+
= \Tc_*(\cpvar{expr}) + \Tc_+(\cpvar{expr})
\end{split}
\end{equation}
if only multiplication and addition are essential for the total time consumed.
For the purpose of analyzing time complexity theoretically instead of estimating practical running time of programs, we can regard the units $\tau_*$ and $\tau_+$ as $1$, thus
the time consumption is equal to the flops involved.

Let $\scc{}{\cpvar{expr}}$ denote the space complexity of computing \cpvar{expr}, 
which means the memories required. For allocating memories for a sequence $\cpvar{v} = \seq{v_0, v_1, \cdots, v_{\ell-1}}$ with positive length $\ell\ge 1$, the space complexity will be $\scc{}{\cpvar{v}} = \ell$. If the memory consumption for the single element of $\cpvar{v}_i$ is $\delta = \cpvar{sizeof}(v_i)$, then the total memories for the sequence of \cpvar{v} will be $\ell\delta$. Just like the analysis of time complexity, the memory unit $\delta$ can be regarded as $1$. Consequently, the key problem of estimating space complexity is to estimate the
counting of memories instead of concrete memory units.  

For the algorithm named with $\ProcName{Alg}$, its TCVC is
\begin{equation}
\begin{split}
\tcc{}{\ProcName{Alg}} 
= \sum_{\cpvar{expr}\in \ProcName{Alg}} \tcc{}{\cpvar{expr}}
=[\tcc{*}{\ProcName{Alg}},\tcc{+}{\ProcName{Alg}}]
\end{split}
\end{equation}
and the time for computation is
\begin{equation}
\begin{split}
\Tc(\ProcName{Alg}) 
= \braket{\tcc{}{\ProcName{Alg}} }{\vec{\tau}} 
= \tcc{*}{\ProcName{Alg}}\tau_*+\tcc{+}{\ProcName{Alg}}\tau_+ 
\end{split}
\end{equation}
Similarly, the space complexity of the  algorithm $\ProcName{Alg}$ is
\begin{equation}
\scc{}{\ProcName{Alg}} 
= \sum_{\cpvar{expr}\in \ProcName{Alg}} \scc{}{\cpvar{expr}}
\end{equation}
If there are some parameters $n, m, \cdots$ for \ProcName{Alg}, then we will take one of the following notations
\begin{equation*}
\tcc{\ProcName{Alg}}{n, m, \cdots},  \Tc_{\ProcName{Alg}}(n, m, \cdots), 
\Tf_{\ProcName{Alg}}(n, m, \cdots),  \scc{\ProcName{Alg}}{n, m, \cdots} 
\end{equation*}
to represent the computational complexity according to practical requirements and interests.

In algorithm analysis, we take the big-O notation \cite{TAOCP-1,CLRS-2022}
\begin{equation}
\BigO{g(n)} 
= \set{h(n): \exists c>0, \exists n_0>0, \forall n \ge n_0, \quad 0\le h(n)\le c\cdot g(n)}
\end{equation}
to represent time or space complexity of interest with an upper bound $g(n)$ where both $c\in \mathbb{R}^+$ and $n_0\in \mathbb{Z}^+$ are constant. In this paper we will use the following fact
$$
\BigO{g(n)\log_2 n} = \BigO{g(n)\mceil{\log_2 n}}, \quad \forall g(n) 
$$
in the analysis of time complexity.

\vspace{2cm}

\noindent\textbf{Citation}: H.-Y. Zhang, Y. Zhou and Z.-Q. Feng. Balanced Binary Tree Schemes for Computing Zernike Radial Polynomials, \textit{IEEE Access}, 11(10): 106567-106579, 2023.  arXiv:2212.02495 [math.NA].  Available online: \blue{\url{https://doi.org/10.1109/ACCESS.2023.3312717}},  \blue{\url{https://arxiv.org/abs/2212.02495}}

\end{document}